\pdfoutput=1

\documentclass[newpagesize]{crscat}

\svnInfo $Id: csimp-intro.tex 345 2009-03-04 12:20:34Z vtl $

\renewcommand{\draftname}{PREPRINT}

\DeclareMathOperator{\MinCtl}{MinCtl}
\DeclareMathOperator{\MaxCtl}{MaxCtl}

\numberwithin{equation}{section}

\title{A simplicial model for proper homotopy types}
\authortrunggwdg

\begin{document}

\begin{abstract}
The singular simplicial set $\Sing(\frakX)$ of a space $\frakX$ completely 
captures its weak homotopy type. We introduce a category of \emph{controlled 
sets}, yielding \emph{simplicial controlled sets}, such that one can 
functorially produce a singular simplicial controlled set 
$\CSing(\MaxCtl(\frakX))$ from a locally compact $\frakX$. We then argue that 
this $\CSing(\MaxCtl(\frakX))$ captures the (weak) \emph{proper} homotopy type 
of $\frakX$. Moreover, our techniques strictly generalize the classical 
simplicial situation: e.g., one obtains, in a unified way, singular homology 
with compact supports and (Borel--Moore) singular homology with locally finite 
supports, as well as the corresponding cohomologies.
\end{abstract}

\maketitle


\section{Introduction}

It is well known that if $\frakX$ is a topological space then its singular 
simplicial set $\Sing(\frakX)$ (which, in dimension $n$, is just the set of 
continuous maps from the topological $n$-simplex to $\frakX$) completely 
captures the weak homotopy type of $\frakX$. The intuitive argument here is 
that $\Sing(\frakX)$ contains all information about maps from spheres (of each 
dimension) to $\frakX$ since spheres can be triangulated, and that it also 
contains all information about homotopies of such maps since homotopies of 
spheres can be realized simplicially (i.e., $S^n \cross \ccitvl{0,1}$ can be 
triangulated). As such, homotopy theory can, if desired, be done in an entirely 
simplicial way.

On the other hand, $\Sing(\frakX)$ contains no information about the proper (or 
locally compact) structure of $\frakX$: the realization $|\Sing(\frakX)|$ of 
$\Sing(\frakX)$ is almost never locally compact, even if $\frakX$ is locally 
compact (or even compact). Thus if one wants to study the (weak) proper 
homotopy invariants (e.g., the Borel--Moore singular homology, various 
analytically obtained invariants) of a locally compact space $\frakX$, it does 
not suffice to only consider the simplicial set $\Sing(\frakX)$.

But of course one can develop, e.g., the locally finitely supported 
Borel--Moore singular homology of a locally compact $\frakX$ in a 
simplicial/singular way, by considering locally finite chains. (The original 
approach of Borel--Moore~\cite{MR0131271} is sheaf-theoretic in nature. The 
approach using locally finite chains is perhaps a bit folkloric, though it has 
been rigorously exposed in \cite{MR1388308}.) Information about whether or not 
a collection of singular simplices is `local finite' is not contained in 
$\Sing(\frakX)$; we solve this problem by replacing the category of sets with a 
category of \emph{controlled sets} (or `sets with supports'), yielding 
\emph{simplicial controlled sets}. In particular, for $\frakX$ locally compact, 
we naturally get a simplicial controlled set $\CSing(\MaxCtl(\frakX))$. A 
feature of our approach is that it simultaneously generalizes the classical 
`compactly supported' case and the proper case: the singular simplicial set 
$\Sing(\frakX)$ can be viewed as a simplicial controlled set 
$\CSing(\MinCtl(\frakX))$ with `minimal control' (see 
Remark~\ref{rmk:top-minctl}), and doing so yields classical results.

After establishing basic definitions and properties, we endeavour to show that 
$\CSing(\MaxCtl(\frakX))$ really does capture the weak proper homotopy type of 
a locally compact $\frakX$. We do not do this in full generality, since doing 
so would require the development of a large amount of theory, in particular a 
corresponding category of `controlled topological spaces' together with closed 
symmetric monoidal structures on these categories (see 
Remark~\ref{rmk:ctl-real-sing}), but only to the extent that is immediately 
possible (but we do give indications as to how to proceed in full generality). 
Just as $\Sing(\frakX)$ contains all information about maps to $\frakX$ from 
finite simplicial complexes (a fact which can be expressed by an adjunction), 
we show that $\CSing(\MaxCtl(\frakX))$ contains all information about 
\emph{proper} maps to $\frakX$ from \emph{locally finite} simplicial complexes.

Finally, we show how the locally finitely supported Borel--Moore singular 
homology (and the corresponding cohomology with compact supports) of $\frakX$ 
can be obtained easily and immediately from $\CSing(\MaxCtl(\frakX))$. Indeed, 
the exact same process applied to $\CSing(\MinCtl(\frakX))$ (which, as 
indicated above, corresponds to the standard singular simplicial set 
$\Sing(\frakX)$) yields the usual singular homology with compact supports (and 
its corresponding cohomology).

We note that our terminology originates in coarse geometry; see in particular 
\cite{MR2007488} (and also the author's~\cite{crscat-I}). It goes without 
saying that there are also applications of our work to coarse geometry, even 
though we do not discuss them here. We have resisted changing our terminology 
of `controlled set' to the arguably more transparent `sets with support' (or 
`supported sets'), since the latter seems to lead to further unwieldy 
terminology.

I would like to thank Bernhard Hanke for his hospitality during my stay at 
Ludwig-Maximilians-Universit\"at M\"unchen (Germany), where most of this paper 
was written. The work was begun at the University of Victoria (Canada), under 
the auspices of Heath Emerson, Marcelo Laca, John Phillips, and Ian Putnam.


\section{Controlled sets}

Denote the category of (small) sets by $\CATSet$. Recall that a set map $f 
\from Y \to X$ is \emph{proper} if $f^{-1}(\set{x})$ is finite for all $x \in 
X$. We have a subcategory $\CATPSet$ of $\CATSet$ with only proper maps, but it 
is rather poorly behaved. We will work with a better-behaved category which 
generalizes both $\CATSet$ and $\CATPSet$.

A \emph{controlled set} $X$ is a set $X$ equipped with a set $\Units(X) 
\subseteq \powerset(X)$ of subsets of $X$ such that:
\begin{enumerate}
\item every finite subset $K$ of $X$ is in $\Units(X)$;
\item if $K$ is in $\Units(X)$, then so are all subsets of $K$; and
\item if $K$ and $K'$ are in $\Units(X)$, then so is $K \union K'$.
\end{enumerate}
The sets in $\Units(X)$ are called the \emph{controlled subsets} of $X$ and 
$\Units(X)$ is called a \emph{control structure} on the set $X$. A set map $f 
\from Y \to X$ between controlled sets is a \emph{controlled map} if:
\begin{enumerate}
\item for each $L$ in $\Units(Y)$, the image $f(L)$ is in $\Units(X)$; and
\item for each $L$ in $\Units(Y)$ and each $x$ in $X$, $f^{-1}(\set{x}) 
    \intersect L$ is finite (i.e., the restriction $f|_L$ is proper).
\end{enumerate}
One can check that (small) controlled sets and controlled maps form a category, 
which we denote by $\CATCSet$.

There is an obvious forgetful functor $\Forget \from \CATCSet \to \CATSet$ 
which has a left adjoint $\MinCtl \from \CATSet \to \CATCSet$ defined putting, 
for any (small) set $S$, $\MinCtl(S) \defeq S$ with $\Units(\MinCtl(S))$ the 
set of all finite subsets of $S$. Not only does $\CATSet$ embed fully and 
faithfully in $\CATCSet$, via the functor $\MinCtl$, but so too does 
$\CATPSet$, via a functor $\MaxCtl \from \CATPSet \to \CATCSet$: for any set 
$S$, put $\MaxCtl(S) \defeq S$ and $\Units(\MaxCtl(S)) \defeq \powerset(S)$.


\section{Simplicial sets and topology}

Denote the (topologists') simplicial category by $\CATSimp$; its objects are 
$\bfzero, \bfone, \bftwo, \dotsc$ (labelled by geometric dimension) and the 
arrows $\bfn \to \bfm$ correspond to the weakly order-preserving maps from the 
total order on $n+1$ elements to the total order on $m+1$ elements. Recall also 
that for each $n$ we have $n+2$ \emph{coface} maps
\[
    d^0, d^1, \dotsc, d^{n+1} \from \bfn \to \bfn \bfplus \bfone
\]
and $n+1$ \emph{codegeneracy} maps
\[
    s^0, s^1, \dotsc, s^n \from \bfn \bfplus \bfone \to \bfn
\]
(we suppress the additional notation required to give the source/targets of 
these maps). Together with the \emph{cosimplicial identities}, these maps 
generate $\CATSimp$.

Recall that a \emph{simplicial object} $X$ in a category $\bfC$ is just an 
object in the functor category $\bfs\bfC \defeq \bfC^{\CATSimp^{\TXTop}}$, 
i.e., a functor $X \from \CATSimp^{\TXTop} \to \bfC$ (or a contravariant 
functor from $\CATSimp$ to $\bfC$). As usual, write $X_n \defeq X(\bfn)$, $d_i 
\defeq X(d^i)$ (the \emph{face} maps of $X$), and $s_j \defeq X(s^j)$ (the 
\emph{degeneracy} maps of $X$). Morphisms in $\bfs\bfC$ are just natural 
transformations, so a map $Y \to X$ between simplicial objects in $\bfC$ is 
given by arrows $Y_n \to X_n$ (in $\bfC$), $n = 0,1,2,\dotsc$.

A \emph{simplicial set} is just a simplicial object in $\CATSet$. For each $n$, 
the \emph{standard $n$-simplex} $\Delta^n$ is the simplicial set represented by 
$\bfn$, i.e., $\Delta^n \defeq \Hom_{\CATSimp}(\blankvar,\bfn)$. The set of 
\emph{$n$-simplices} of a simplicial set $X$ is $X_n \defeq X(\bfn)$, or 
equivalently (by the Yoneda Lemma) $X_n \cong \Hom_{\CATsSet}(\Delta^n,X) = 
\Nat(\Delta^n,X)$, where as usual $\CATsSet \defeq 
\CATSet^{\CATSimp^{\TXTop}}$; a \emph{simplex} of $X$ is an element $x \in X_n$ 
for some $n$ (strictly speaking, one must keep track of the dimension 
explicitly since the $X_n$ need not be disjoint).

The \emph{realization functor} is a functor $|\blankvar| \from \CATSimp \to 
\CATTop$, where $\CATTop$ is some `convenient category' of topological spaces 
such as the category of compactly generated spaces (see, e.g., 
\cite{MR1702278}*{Chapter~5}). On objects, it is given by
\[
    |\bfn| \defeq \set{(x_0,\dotsc,x_n) \in (\setRplus)^{n+1}
            \suchthat x_0 + \dotsb + x_n = 1} \subseteq \setR^{n+1}.
\]
If $X$ is a simplicial set, we obtain its \emph{realization} $|X|$ by the coend 
(in $\CATTop$)
\begin{equation}\label{eq:real}
    |X| \defeq \int^{\bfn} \Top(X(\bfn)) \tensor |\bfn|
\end{equation}
where $\Top \from \CATSet \to \CATTop$ is the `discrete space' functor and 
$\tensor$ is just the ordinary cartesian product (since $\CATTop$ is cartesian 
closed). This can be made into a functor $|\blankvar| \from \CATsSet \to 
\CATTop$ (again called the \emph{realization functor}). There are a number of 
ways of describing the realization $|X|$ but they are all equivalent and all 
yield Milnor's geometric realization~\cite{MR0084138}.

Conversely, given a topological space $\frakX$, we can form its \emph{singular 
simplicial set}
\[
    \Sing(\frakX) \defeq \Hom_{\CATTop}(|\blankvar|,\frakX).
\]
It is clear that this yields a functor $\Sing \from \CATTop \to \CATsSet$. The 
basic, but essential, result is that realization is left adjoint to $\Sing$, 
i.e., for all spaces $\frakX$ and simplicial sets $Y$,
\begin{equation}\label{eq:real-sing}
    \Hom_{\CATTop}(|Y|,\frakX) \cong \Hom_{\CATsSet}(Y,\Sing(\frakX))
\end{equation}
(isomorphism of sets, natural in both $\frakX$ and $Y$).

Recall that a simplex of a simplicial set is \emph{nondegenerate} if it is not 
in the image of any degeneracy map. A simplicial set $X$ is \emph{finite} if it 
only has finitely many nondegenerate simplices. If $X$ is finite, then each 
$X_n$ must be a finite set.

\begin{remark}
The definition of `finite simplicial set' is a bit delicate: one cannot just 
demand that $X$ only have finitely many simplices, since this can only hold if 
$X_n = \emptyset$ for all $n$. On the other hand, it is not enough to insist 
that each $X_n$ be finite, since $X$ could then still be infinite dimensional.
\end{remark}

Since the realization $|X|$ is a $CW$-complex with one cell for each 
nondegenerate simplex~\cite{MR0084138}*{Theorem~1}, the following is clear and 
well known.

\begin{proposition}\label{prop:fin-simp-cpt}
If $X$ is a finite simplicial set, then its realization $|X|$ is compact, 
Hausdorff.
\end{proposition}

Simplices $x$, $x'$ of a simplicial set $X$ are \emph{adjacent} if there exist 
arrows $f$, $f'$ in $\CATSimp$ such that $X(f)(x) = X(f')(x')$. A simplicial 
set $X$ is \emph{locally finite} if each simplex of $X$ is adjacent to only 
finitely many nondegenerate simplices. If $X$ is locally finite, then all the 
maps with which $X$ comes are proper (this follows easily from the injectivity  
of the degeneracy maps), i.e., $X$ is actually a simplicial object in 
$\CATPSet$.

\begin{remark}
It is not sufficient for our purposes to only insist that $X$ be a simplicial 
object in $\CATPSet$ (but note that others, e.g., \cite{MR1388308}, define 
`locally finite simplicial set' as such). Indeed, our definition essentially 
includes an assertion of `local finite dimensionality'. Our definition actually 
reduces to the a-priori-weaker requirement that each $0$-simplex be adjacent to 
only finitely many nondegenerate simplices, since if two simplices are adjacent 
then they are both adjacent to some common $0$-simplex. Observe also that if 
$X$ is locally finite and finite in each dimension, then it is actually finite.
\end{remark}

We have the following analogue of Proposition~\ref{prop:fin-simp-cpt}, which 
should also be clear.

\begin{proposition}\label{prop:lf-simp-lc}
If $X$ is a locally finite simplicial set, then its realization $|X|$ is 
locally compact, Hausdorff.
\end{proposition}

In fact, more is true. The category $\CATsPSet \defeq 
\CATPSet^{\CATSimp^{\TXTop}}$ of simplicial objects in $\CATPSet$ maps 
naturally to the category $\CATsSet$ of simplicial sets, and we have already 
noted that locally finite simplicial sets `lift' to $\CATsPSet$. Say that a 
morphism between locally finite simplicial sets is \emph{proper} if it lifts to 
a morphism in $\CATsPSet$ (i.e., if all the set maps involved are proper). We 
then get the following.

\begin{proposition}\label{prop:prop-prop}
If $f \from Y \to X$ is a proper morphism between locally finite simplicial 
sets, then the induced map $|f| \from |Y| \to |X|$ on realizations is proper.
\end{proposition}

Denote the category of locally compact, Hausdorff spaces and proper maps by 
$\CATPTop$ (a subcategory of $\CATTop$). Combining 
Propositions~\ref{prop:lf-simp-lc} and~\ref{prop:prop-prop}, we get the 
following.

\begin{theorem}
The usual realization functor yields a functor (which we continue to denote by 
$|\blankvar|$) from the category of locally finite simplicial sets and proper 
simplicial maps (which can be considered as a full subcategory of $\CATsPSet$) 
to $\CATPTop$.
\end{theorem}


\section{Simplicial controlled sets and topology}

Our main goal is to generalize \eqref{eq:real-sing} to the controlled/`proper' 
setting. The essential difficulty in doing so is that $|\Sing(\frakX)|$ is 
almost never locally compact, so the canonical maps $\epsilon_{\frakX} \from 
|\Sing(\frakX)| \to \frakX$ (given by the counit of adjunction) in 
\eqref{eq:real-sing} are almost never proper. In this paper, we will only prove 
a special case of the controlled version of \eqref{eq:real-sing}, since giving 
a general adjunction would take us too far afield.

\begin{definition}
A \emph{simplicial controlled set} is a simplicial object in $\CATCSet$ (i.e., 
an object of $\CATsCSet \defeq \CATCSet^{\CATSimp^{\TXTop}}$), and a morphism 
between simplicial controlled sets is just a morphism in $\CATsCSet$ (i.e., a 
natural transformation).
\end{definition}

Define a functor $\CSing \circ \MaxCtl \from \CATPTop \to \CATsCSet$ (our 
notation is justified in Remark~\ref{rmk:ctl-real-sing}) as follows. The 
forgetful functor $\Forget \from \CATCSet \to \CATSet$ naturally yields a 
functor $\Forget_{\#}(\CSing \circ \MaxCtl) \from \CATPTop \to \CATsSet$, which 
we insist should coincide with our original functor $\Sing$ (composed with the 
inclusion $\CATPTop \injto \CATTop$). That is, if $\frakX$ is a locally 
compact, the sets underlying the simplicial controlled set 
$\CSing(\MaxCtl(\frakX))$ are just the sets of the simplicial set 
$\Sing(\frakX)$, and similarly for the underlying set maps. We must also give 
the control structures $\Units(\CSing(\MaxCtl(\frakX))_n)$: say that $S 
\subseteq \CSing(\MaxCtl(\frakX))_n = \Sing(\frakX)_n$ is a controlled subset 
if it is \emph{locally finite} in the sense that each point of $\frakX$ has a 
neighbourhood which meets (the images of) only finitely many simplices in $S$. 
It is easy to check that this does indeed define a functor (in particular, that 
all maps are controlled).

For any simplicial set $X$, $\MinCtl \circ X$ is a simplicial controlled set. 
(Our compositional notation $\MinCtl \circ X$ is due to the fact that $X$ is a 
functor. Concretely, $(\MinCtl \circ X)_n = \MinCtl(X_n)$.) This embeds the 
category of simplicial sets fully and faithfully in the category of simplicial 
controlled sets. By considering locally finite simplicial sets to be simplicial 
objects in $\CATPSet$, we can likewise embed the category of locally finite 
simplicial sets and proper morphisms in $\CATsCSet$: send each locally finite 
simplicial set $X$ to the simplicial controlled set $\MaxCtl \circ X$.

Our main theorem is the following.

\begin{theorem}\label{thm:prop-real-sing}
For each locally compact, Hausdorff space $\frakX$ and each locally finite 
simplicial set $Y$, we have an isomorphism
\begin{equation}\label{eq:prop-real-sing}
    \Hom_{\CATPTop}(|Y|,\frakX)
        \cong \Hom_{\CATsCSet}(\MaxCtl \circ Y, \CSing(\MaxCtl(\frakX)))
\end{equation}
(where, as before, $\CATPTop$ is the category of locally compact, Hausdorff 
spaces and proper maps). These isomorphisms are natural in $\frakX$ (and proper 
maps) and $Y$ (and proper maps), and the isomorphisms are compatible with the 
classical adjunction \eqref{eq:real-sing}.
\end{theorem}

More precisely, `compatible' means that the square
\[\begin{CD}
    \Hom_{\CATPTop}(|Y|,\frakX)
        @<{\sim}-> \Hom_{\CATsCSet}(\MaxCtl \circ Y, \CSing(\MaxCtl(\frakX)))
        \\
    @VVV @VVV \\
    \Hom_{\CATTop}(|Y|,\frakX) @<{\sim}-> \Hom_{\CATsSet}(Y,\Sing(\frakX))
\end{CD}\]
commutes (above, we omitted some inclusions of categories).

\begin{remark}\label{rmk:ctl-real-sing}
The isomorphisms \eqref{eq:prop-real-sing} do not constitute an adjunction 
since $\CSing \circ \MaxCtl$ does not produce locally finite simplicial sets. 
The complete, correct way to proceed is as follows:
\begin{enumerate}
\item Make $\CATCSet$ into a closed (symmetric monoidal) category, i.e., equip 
    it with a `tensor product' operation $\tensor$ (in the obvious way) to make 
    it into a symmetric monoidal category and then construct a compatible 
    internal $\Hom$.
\item Define a category $\CATCTop$ of \emph{controlled topological spaces} 
    along with a `forgetful' functor $\Forget \from \CATCTop \to \CATCSet$, and 
    make it into a closed category in a way compatible with the structures on 
    $\CATCSet$. Essentially, the objects of $\CATCTop$ are spaces equipped with 
    a suitable family of \emph{controlled subspaces} (or `supports') and the 
    maps are continuous and proper relative to the controlled subspaces.
\item Define various functors compatible with the existing ones:
    \[
    \begin{CD}
        \CATSet @>{\Top}>> \CATTop \\
        @V{\MinCtl}VV @VV{\MinCtl}V \\
        \CATCSet @>{\CTop}>> \CATCTop
    \end{CD}
    \qquad\quad\text{and}\quad\qquad
    \begin{CD}
        \CATPSet @>{\PTop}>> \CATPTop \\
        @V{\MaxCtl}VV @VV{\MaxCtl}V \\
        \CATCSet @>{\CTop}>> \CATCTop
    \end{CD}
    \]
    (where $\PTop \from \CATPSet \to \CATPTop$ is the obvious functor).
\item Define $\CSing \from \CATCTop \to \CATsCSet$ by putting, for $\frakX$ a 
    controlled space,
    \[
        \CSing(\frakX) \defeq \Forget(\Hom_{\CATCTop}(|\blankvar|,\frakX))
    \]
    (where $\Hom_{\CATCTop}$ is actually the internal $\Hom$ of $\CATCTop$, 
    whence we get a controlled set by `forgetting'). One should check that this 
    agrees with our explicit construction of $\CSing \circ \MaxCtl$ above and 
    that it generalizes $\Sing$ in the sense that $\Sing = \Forget_{\#}(\CSing 
    \circ \MinCtl)$ (notation as above; see also Remark~\ref{rmk:top-minctl} 
    below).
\item Define a realization functor $|\blankvar| \from \CATsCSet \to \CATCTop$ 
    by the coend (in $\CATCTop$)
    \[
        |X| \defeq \int^{\bfn} \CTop(X(\bfn)) \tensor |\bfn|
    \]
    for each simplicial controlled set $X$ (strictly speaking, $|\bfn|$ should 
    be written as $\MinCtl(|\bfn|)$; see Remark~\ref{rmk:top-minctl} below). 
    This is entirely analogous to the definition \eqref{eq:real} of the 
    classical realization of a simplicial set, but one should verify that it 
    generalizes the classical realization (in the sense that 
    $\Forget(|\MinCtl_{\#}(\blankvar)|)$ is just the classical realization 
    $|\blankvar| \from \CATsSet \to \CATTop$).
\item Prove that realization is left adjoint to $\CSing$.
\end{enumerate}
Unfortunately, the second step, constructing $\CATCTop$ involves nontrivial 
work in general topology, with many of the same difficulties entailed in 
finding a `convenient category' of topological spaces (as well as other 
difficulties). We thus leave this task to the future.
\end{remark}

\begin{remark}\label{rmk:top-minctl}
By analogy with $\CSing \circ \MaxCtl$, we may also define
\[
    \CSing \circ \MinCtl \defeq \MinCtl_{\#}(\Sing) \from \CATTop \to \CATCSet,
\]
i.e., for any space $\frakX$, $\CSing(\MinCtl(\frakX)) \defeq \MinCtl \circ 
\Sing(\frakX)$. Here, the notation indicates that to any space $\frakX$ can be 
made into a controlled topological space $\MinCtl(\frakX)$ by equipping it with 
its family of compact subspaces.
\end{remark}

\begin{proof}[Proof of Theorem~\ref{thm:prop-real-sing}]
Fix $\frakX$ and $Y$. Since we insisted on compatibility with the classical 
version, the adjunction \eqref{eq:real-sing} actually determines the 
isomorphism \eqref{eq:prop-real-sing} and ensures naturality; it is only a 
matter of checking that the maps which \eqref{eq:real-sing} provides are proper 
or controlled (as appropriate).

Suppose that $f \from |Y| \to \frakX$ is proper. We get a natural map $F = 
(F_n) \from Y \to \Sing(\frakX)$ of simplicial sets, i.e., for each $n$ we get 
a set map $F_n \from Y_n \to \Sing(\frakX)_n$. We must check that the $F_n$ are 
controlled with respect to the control structures of $\MaxCtl(Y_n)$ and 
$\CSing(\MaxCtl(\frakX))$, i.e., that the $F_n$ are proper and have controlled 
images. Since $f$ is proper, for each compact $K \subseteq \frakX$, $f^{-1}(K)$ 
intersects (the realizations of) only finitely many nondegenerate simplices; 
thus each compact $K$ intersects the images of only finitely many nondegenerate 
simplices of $F(Y)$. It follows that in each dimension $n$, each $K$ intersects 
the images of only finitely many simplices, so each $F_n(Y_n)$ is controlled. 
The same kind of argument also shows that each $F_n$ is proper.

Conversely, suppose that $F = (F_n) \from \MaxCtl \circ Y \to 
\CSing(\MaxCtl(\frakX))$ is a map of simplicial controlled sets (hence a map of 
simplicial sets). Naturally, we get a continuous map $f \from |Y| \to \frakX$. 
We must check that $f$ is proper, so fix $K \subseteq \frakX$ compact. By 
properness of the $F_n$, $f^{-1}(K)$ can intersect the (realizations of) 
finitely many $n$-simplices of $Y$. Then, by local finiteness of $Y$, 
$f^{-1}(K)$ can only intersect finitely many nondegenerate simplices.
\end{proof}


\section{Homology and cohomology}

Let us briefly indicate how we can develop homology and cohomology in our 
setting. Fix an abelian group $G$ (of coefficients), e.g., $G \defeq \setZ$. 
Given a controlled set $X$, let $G[X]$ be the set (indeed, abelian group) of 
formal sums
\[
    \sum_{x \in X} [x] g_x, g_x \in G,
\]
with controlled support, i.e., $\set{x \in X \suchthat g_x \neq 0}$ is in 
$\Units(X)$. We make $G[\blankvar]$ into \emph{covariant} functor (from 
$\CATCSet$ to the category $\CATAb$ of abelian groups) in the way suggested by 
the summation notation: if $f \from Y \to X$ is a controlled map between 
controlled sets, then
\[
    f_*\bigl({\textstyle\sum [y] g_y}\bigr) \defeq {\textstyle\sum [f(y)] g_y}
\]
(adding coefficients when necessary). Formally, $G[X]$ is just a set of 
functions $X \to G$ (or, equivalently, a subgroup of the product $\prod_{x \in 
X} G$), but this is misleading in terms of functoriality.

Thus, if $X$ is a \emph{simplicial controlled set}, then $G[X]$ (strictly 
speaking, $G[\blankvar] \circ X$) is a simplicial abelian group. We then get a 
chain complex and homology groups $H_{\grstar}(G[X])$ in the usual way. We 
immediately get the following (see, e.g., \cite{MR1388308}).

\begin{proposition}
For any locally compact space $\frakX$, $H_{\grstar}(G[X])$ is just the 
(locally finitely supported) Borel--Moore singular homology of $\frakX$.
\end{proposition}

Even more trivially, we can recover the usual singular homology.

\begin{proposition}
For any space $\frakX$, $H_{\grstar}(G[\CSing(\MinCtl(\frakX))])$ (notation as 
in Remark~\ref{rmk:top-minctl}) is just the usual (compactly supported) 
singular homology of $\frakX$.
\end{proposition}

To get cohomology, we give a contravariant functor from $\CATCSet$ to $\CATAb$ 
(i.e., a covariant functor $\CATCSet^{\TXTop} \to \CATAb$). Say that a subset 
$S$ of a controlled set $X$ is \emph{cocontrolled} if $S \intersect T$ is 
finite for all $T \in \Units(X)$. Let $G[X]^*$ be the set (indeed, abelian 
group) of functions $X \to G$ with cocontrolled support. One can then check 
that this yields a \emph{contravariant} functor (by pulling back functions, as 
usual).

For any simplicial controlled set $X$, $G[X]^*$ (strictly speaking, 
$G[\blankvar]^* \circ X$) is a cosimplicial abelian group, and we get a cochain 
complex and cohomology groups $H^{\grstar}(G[X]^*)$. One can check the 
following results for cohomology.

\begin{proposition}
For any locally compact space $\frakX$, 
$H^{\grstar}(G[\CSing(\MaxCtl(\frakX))]^*)$ is the usual singular homology of 
$\frakX$ with compact supports.
\end{proposition}

\begin{proposition}
For any space $\frakX$, $H^{\grstar}(G[\CSing(\MinCtl(\frakX))]^*)$ is the 
usual singular cohomology of $\frakX$.
\end{proposition}

\begin{remark}
For any controlled set $X$ and abelian groups $G$ and $H$, observe that there 
is an obvious pairing
\[
    (\Hom_{\CATAb}(H,G))[X]^* \cross H[X] \to G,
\]
namely the one given by
\[
    \bigl\langle \alpha, \sigma \bigr\rangle
        \defeq \sum_x \alpha(x)(h_x)
\]
for $\alpha \in (\Hom(H,G))[X]^*$ (hence $\alpha$ is a function $X \to 
\Hom(H,G)$) and $\sigma \defeq \sum [x] h_x \in H[X]$; the sum is actually 
finite since $\alpha$ has cocontrolled support and $\sigma$ has controlled 
support. This evidently induces the usual pairing between homology and 
cohomology.
\end{remark}

\begin{remark}
One may obtain the contravariant functor $G[\blankvar]^*$ in a slightly more 
satisfactory way as follows. First, observe that, for any $X$, one can give 
$G[X]$ slightly more structure than just its structure as an abelian group. 
Given a collection (possibly infinite) $\set{\sigma_i}_{i \in I}$ ($I$ some 
index set) of elements of $G[X]$, one can form the (possibly infinite) sum 
$\sum_{i \in I} \sigma_i$ if each $x \in X$ is contained in only finitely many 
of the supports of the $\sigma_i$. That is, $G[X]$ has an additional structure 
which allows one to take certain infinite sums.

One can axiomatize such objects, namely abelian groups in which certain 
infinite sums are defined, and obtain a category $\CATCAb$ of \emph{controlled 
abelian groups}. Then one can make $\CATCAb$ into a closed category, in 
particular giving it an internal $\Hom$. Giving $G$ the structure of a 
controlled abelian group with `minimal control' (i.e., only finite sums are 
defined), we get a controlled abelian group $\Hom_{\CATCAb}(\setZ[X],G)$ which 
turns out to be our explicitly-constructed exactly $G[X]^*$.
\end{remark}



\end{document}